\documentclass[12pt]{amsart}
\usepackage{amssymb}

\newtheorem{thm}{Theorem}[section] 
\newtheorem{prp}{Proposition}[section] 
\newtheorem{lem}{Lemma}[section] 
\newtheorem{cor}{Corollary}[section] 
\begin{document}
\author{G\'erard Endimioni}
\address{C.M.I-Universit\'{e} de Provence, 
39, rue F. Joliot-Curie, F-13453 Marseille Cedex 13}
\email{endimion@gyptis.univ-mrs.fr}
\title[Normal automorphisms]{Normal automorphisms of a free metabelian nilpotent group}
\subjclass[2000]{20E36, 20F28.}
\keywords{Normal Automorphism; Free Metabelian Nilpotent Group}
\begin{abstract} 
An automorphism $\varphi$ of a group $G$ is said to be normal
if $\varphi(H)=H$ for each normal subgroup $H$ of $G$.
These automorphisms form a group containing the group of inner 
automorphisms.
When $G$ is a nonabelian free (or free soluble) group, it is known that
these groups of automorphisms coincide,
but this is not always true when 
$G$ is a free metabelian nilpotent group. The aim of this paper is to 
determine the group of normal automorphisms in this last case.
\end{abstract}
\maketitle
%
\section{Preliminary results}
\noindent In a group $G$, consider a map $\varphi: G\to G$ of the form
$$\varphi:x\to 
x\left[x,u_{1}\right]^{\lambda(1)}\ldots\left[x,u_{m}\right]^{\lambda(m)},$$
where $u_{1},\ldots,u_{m}$ are elements 
of $G$, the exponents $\lambda(1),\ldots,\lambda(m)$ being integers (as usual, the 
commutator $[a,b]$ is defined by $[a,b]=a^{-1}b^{-1}ab$).
When $G$ is metabelian, using the relation 
$[xy,u]=y^{-1}[x,u]y[y,u]$, it is easy to see that $\varphi$ is an 
endomorphism of group. These endomorphisms appear in \cite{KU} as a 
solution to a problem of extension from $G'$ 
to $G$ for certain endomorphisms (also see \cite{CS}).

Now consider in an arbitrary group $G$ the map
$\varphi_{u,v}:x\to x[x,u][x,v]$, with $u,v\in G$ and suppose that $\varphi_{u,v}$ is 
an endomorphism. It is then easy to show that the equality
$\varphi_{u,v}(xy^{-1})=\varphi_{u,v}(x)\varphi_{u,v}(y^{-1})$ is equivalent to the 
relation $[u,y][x,v]=[x,v][u,y]$. Consequently, if $\varphi_{u,v}$ is 
an endomorphism for any $u,v\in G$, the group $G$ is metabelian.
One can summarize these preliminary remarks like this.
\begin{prp}
\bigskip 
In a group $G$, the following 
conditions are equivalent:\\
{\rm (i)} every map $\varphi: G\to G$ of the form
$\displaystyle\varphi:x\to 
x\prod_{i=1}^m\left[x,u_{i}\right]^{\lambda(i)}$
$(u_{i}\in G,\: \lambda(i)\in {\mathbb Z})$
is an endomorphism;\\
{\rm (ii)} $G$ is metabelian.
\end{prp}
It is very easy to see that in an metabelian group, such endomorphisms 
are not necessarily 
automorphisms (consider for example the endomorphism
$x\to x[x,b]^{-1}$ in the nonabelian group of order 6 
defined by the presentation
$\langle a,b\, \vert\, a^3=b^2=(ab)^2=1\rangle$). 
But in a nilpotent group, each map of 
the form 
$$x\to w_{0}x^{\lambda(1)}w_{1}x^{\lambda(2)}\ldots 
x^{\lambda(n)}w_{n}\;\;\; ({\rm with}\; 
\lambda(1)+\lambda(2)+\cdots+\lambda(n)=\pm 1)$$ 
is bijective \cite[Theorem 1]{EN1}. Hence we have:
\begin{prp}
In a metabelian nilpotent group $G$, every map $\varphi: G\to G$ of the form
$\displaystyle\varphi:x\to 
x\prod_{i=1}^m\left[x,u_{i}\right]^{\lambda(i)}$
$(u_{i}\in G,\: \lambda(i)\in {\mathbb Z})$
is an automorphism.
\end{prp}
For convenience sake, in a metabelian nilpotent group, an 
automorphism of the form $\displaystyle x\to 
x\prod_{i=1}^m\left[x,u_{i}\right]^{\lambda(i)}$
will be called a {\em generalized inner automorphism}.
Obviously, every inner automorphism 
is a generalized inner automorphism, since $u^{-1}xu=x[x,u]$.
Notice that in a nilpotent group $G$ of class $\leq 2$, 
we may write 
$$\varphi(x)=x\prod_{i=1}^m\left[x,u_{i}\right]^{\lambda(i)}
=x\prod_{i=1}^m\left[x,u_{i}^{\lambda(i)}\right]
=x\left[x,\prod_{i=1}^m u_{i}^{\lambda(i)}\right].$$
Hence $\varphi(x)=x[x,u]=u^{-1}xu$, where 
$\displaystyle u=\prod_{i=1}^m u_{i}^{\lambda(i)}$.
Consequently, in a nilpotent group of class $\leq 2$, the notions 
of inner automorphism and generalized inner automorphism coincide. 
On the other hand, in the free nilpotent group of class 3 and rank 2 
freely generated by $a$ and $b$, it is easy to see that the map 
$x\to x[x,a,a]$
is a generalized inner automorphism but is not an inner automorphism.

As usual, in a group, the left-normed commutator 
$[x_{1},\ldots,x_{n}]$ is defined inductively by
\begin{eqnarray*}
	[x_{1},\ldots,x_{n}] & = & [x_{1},\ldots,x_{n-1}]^{-1}[x_{1},\ldots,x_{n-1}]^{x_{n}} \\
	 & = & [x_{1},\ldots,x_{n-1}]^{-1}x_{n}^{-1}[x_{1},\ldots,x_{n-1}]x_{n}.
\end{eqnarray*}
The next technical result will be useful in the following.
\begin{prp}
In a group $G$, consider a map
$\varphi:G\to G$ of the form
$$\varphi(x)=x\prod_{i=1}^{n}[x,v_{i,1},\ldots,v_{i,\sigma(i)}]^{\eta(i)} 
\quad ({\eta(i)}\in {\mathbb Z}),$$
for some function $\sigma :\{ 1,\ldots,n\}\to {\mathbb N}\setminus\{ 0\}$ and 
elements $v_{i,j}\in G$ ($1\leq i\leq n$, $1\leq j\leq \sigma(i)$).
Then $\varphi(x)$ can be written in the form
$$\varphi(x)=
x\left[x,u_{1}\right]^{\lambda(1)}\ldots\left[x,u_{m}\right]^{\lambda(m)}
\quad (\lambda(i)\in {\mathbb Z},\: u_{i}\in G).$$
\end{prp}
\begin{proof} By induction, using the relation 
$[x,y,z]=[x,y]^{-1}[x,z]^{-1}[x,yz]$
\end{proof} 
Frequently in this paper we shall make use of well-known 
commutator identities (see for example \cite[5.1.5]{RO}). In particular,
we have the following relations, 
valid in a metabelian group $G$, for any $x,y,z\in G$, $t\in G'$ and 
$\lambda\in {\mathbb Z}$:
\begin{center}
\begin{tabular}{cccc}
	$[xt,y]=[x,y][t,y]$ & {} & {} & $[t^{\lambda},y]=[t,y]^{\lambda}$  \\
	$[x,y,z][y,z,x][z,x,y]=1$ & {}  & {} & $[t,x,y]=[t,y,x]$  \\
\end{tabular}
\end{center}
\begin{prp}
The set of generalized inner automorphisms of a 
metabelian nilpotent group $G$ forms a 
(normal) subgroup of the group of automorphisms of $G$.
\end{prp}
\begin{proof} If $\varphi$ and $\psi$ are generalized inner 
automorphisms respectively defined by 
$$\varphi(x)=x\prod_{i=1}^m [x,u_{i}]^{\lambda(i)},\;
\psi(x)=x\prod_{i=1}^n [x,v_{i}]^{\mu(i)},$$
an easy calculus shows that
$$\psi\circ\varphi(x)=x\prod_{i=1}^m [x,u_{i}]^{\lambda(i)}
\prod_{i=1}^n [x,v_{i}]^{\mu(i)}
\prod_{i=1}^m \prod_{j=1}^n [x,u_{i},v_{j}]^{\lambda(i)\mu(j)}.$$
Thus $\psi\circ\varphi$ is a generalized inner automorphism by
Proposition 1.3.\\
It remains to prove that $\varphi^{-1}$ is a 
generalized inner automorphism. 
For that, it suffices to construct for each integer 
$k\geq 1$ a generalized inner automorphism $\psi_{k}$ such that
$\psi_{k}\circ \varphi$ is of the form
$$\psi_{k}\circ \varphi :x\to 
x\prod_{i=1}^{m}[x,v_{i,1},\ldots,v_{i,\sigma(i)}]^{\eta(i)}$$
for some function $\sigma :\{ 1,\ldots,m\}\to {\mathbb N}\setminus\{ 0\}$ and 
elements $v_{i,j}\in G$ ($1\leq i\leq m$, $1\leq j\leq \sigma(i)$),
and where each commutator is of weight $\geq 1+2^{k-1}$ 
(namely, $\sigma(i)\geq 2^{k-1}$ for $i=1,\ldots,m$).
Indeed, since $G$ is nilpotent, this implies that 
$\psi_{k}\circ \varphi(x)=x$ for $k$ large enough, thus 
$\varphi^{-1}=\psi_{k}$ is a generalized inner automorphism, as required.
We argue by induction on $k$. The result is clear when $k=1$
by taking for $\psi_{1}$ the identity map. Now
suppose that for some integer $k\geq 1$, there exists 
a generalized inner automorphism $\psi_{k}$ such that
$\displaystyle \psi_{k}\circ \varphi(x)= 
x\prod_{i=1}^{m}[x,v_{i,1},\ldots,v_{i,\sigma(i)}]^{\eta(i)}$, with
$\sigma(i)\geq 2^{k-1}$ for $i=1,\ldots,m$. Put $\psi_{k+1}=\psi'\circ\psi_{k}$, 
where $\psi'$ is defined by
$\displaystyle \psi'(x)= 
x\prod_{i=1}^{m}[x,v_{i,1},\ldots,v_{i,\sigma(i)}]^{-\eta(i)}$.
We have
$$\psi_{k+1}\circ\varphi(x)=x\prod_{i=1}^{m}[x,v_{i,1},\ldots,v_{i,\sigma(i)}]^{\eta(i)}$$
$$\times \prod_{j=1}^{m}\left[x\prod_{i=1}^{m}[x,v_{i,1},\ldots,v_{i,\sigma(i)}]^{\eta(i)},
v_{j,1},\ldots,v_{j,\sigma(j)}\right]^{-\eta(j)}.$$
Since
$$	\prod_{j=1}^{m}\left[x\prod_{i=1}^{m}[x,v_{i,1},\ldots,v_{i,\sigma(i)}]^{\eta(i)},
v_{j,1},\ldots,v_{j,\sigma(j)}\right]^{-\eta(j)}$$
$$=\prod_{i=1}^{m}[x,v_{j,1},\ldots,v_{j,\sigma(j)}]^{-\eta(j)}
\prod_{j=1}^{m}\prod_{i=1}^{m}[x,v_{i,1},\ldots,v_{i,\sigma(i)},
v_{j,1},\ldots,v_{j,\sigma(j)}]^{-{\eta(i)}\eta(j)}, $$
we obtain
$$\psi_{k+1}\circ\varphi(x)=x\prod_{j=1}^{m}\prod_{i=1}^{m}[x,v_{i,1},\ldots,v_{i,\sigma(i)},
v_{j,1},\ldots,v_{j,\sigma(j)}]^{-{\eta(i)}\eta(j)}$$
and this completes the proof of the proposition.
\end{proof}
%
\section{Main result}
We recall that a normal automorphism $\varphi$ of a group $G$ is an 
automorphism such that $\varphi(H)=H$ for each normal subgroup 
$H$ of $G$. These automorphisms form a subgroup of the group of all 
automorphisms of $G$. Obviously, this subgroup contains the subgroup 
of inner automorphisms of $G$. It happens these subgroups coincide,
for instance when $G$ is a nonabelian free group \cite{LU}, 
a nonabelian free  soluble group \cite{ROM}, 
or a nonabelian free nilpotent group of class $2$ \cite{EN2}. 
On the other hand, the subgroup of inner automorphisms 
is of infinite index in the group of normal automorphisms when $G$ 
is a nonabelian free nilpotent group of class $k\geq 3$ \cite{EN2}.
Also note there are exactly two normal
automorphisms in a (nontrivial) free abelian group: $x\to x$ and 
$x\to x^{-1}$.
The determination of all normal automorphisms 
in a free nilpotent group of class $k\geq 3$ seems to be an open problem.
Note however that for $k=3$, Theorem 2.1 below provides such a 
determination, since the free nilpotent metabelian group of class $3$  
(and of given rank) coincide with the corresponding free nilpotent group of class 
$3$.
One can find in \cite{FG} information about the structure of the group
of normal automorphisms of a nilpotent group. 

Certainly, in a metabelian nilpotent group, each
generalized inner automorphism is a normal automorphism, but 
a normal automorphism need not to be a generalized inner automorphism.  
However, our main result states that  the converse holds in a 
nonabelian free metabelian nilpotent group.
\begin{thm}
In a nonabelian free metabelian nilpotent 
group, the group of normal automorphisms coincides with the group of 
generalized inner automorphisms.
\end{thm}
Before to tackle the proof of this theorem, we give a few consequences.
Let $\varphi$ be an automorphism of a group $G$. According to Schweigert \cite{SC}, 
one says that $\varphi$ is a polynomial automorphism of $G$
if there exist integers $\epsilon_{1},\ldots,\epsilon_{m}\in{\mathbb Z}$ 
and elements $u_{1},\ldots,u_{m}\in G$ such that
$$\varphi(x)=(u_{1}^{-1}x^{\epsilon_{1}}u_{1})\ldots
(u_{m}^{-1}x^{\epsilon_{m}}u_{m})$$
for all $x\in G$. For instance, in a metabelian nilpotent group, any 
generalized inner automorphism is a polynomial automorphism. It turns out that in a
nonabelian free metabelian nilpotent group, these notions coincide.
\begin{cor} An automorphism $\varphi$ of a nonabelian free 
metabelian nilpotent group $G$ is a generalized inner automorphism if and only if 
it is polynomial.
\end{cor}
\begin{proof} Only the part 'if' needs a proof. Therefore, suppose that $\varphi:G\to G$ is a polynomial automorphism
defined by  
$$\varphi(x)=(u_{1}^{-1}x^{\epsilon_{1}}u_{1})\ldots
(u_{m}^{-1}x^{\epsilon_{m}}u_{m})$$ 
and put 
$\epsilon= \epsilon_{1}+\cdots+\epsilon_{m}$.
If $\gamma_2(G)=[G',G]$ is the second term of the lower central series of $G$, 
it is not difficult to see that $\varphi$ induces the automorphism $\overline{\varphi}:x\mapsto x^{\epsilon}[x,u]$ 
on the free nilpotent group $G/\gamma_2(G)$, where
$u=u_1^{\epsilon_{1}}\ldots u_m^{\epsilon_{m}}\gamma_2(G)$. An easy calculation
shows that the relation $\overline{\varphi}(xy)=\overline{\varphi}(x)\overline{\varphi}(y)$
is equivalent to the relation $(xy)^{\epsilon}=x^{\epsilon}y^{\epsilon}$. This implies that $\epsilon=1$, and so $\varphi^{-1}$ is polynomial by
\cite[Theorem 1]{EN1}. Now it is clear that $\varphi$ is a normal automorphism, hence
$\varphi$ is a generalized inner automorphism by Theorem 2.1. 
\end{proof}
In \cite{EN3}, it is proved that the group generated by all  
polynomial automorphisms of a metabelian group is itself metabelian.
Therefore, we can state:
\begin{cor} The group of normal automorphisms of a free 
metabelian nilpotent group is metabelian.
\end{cor}
Note that this result is valid even if the relatively free group is abelian, since the
group of normal automorphisms is abelian in this case.\\ 
Recall that an automorphism of a group $G$ is said to be an IA-automorphism if it induces the identity automorphism on $G/G'$. We shall write $\textrm{IA}(G)$ for the group of
IA-automorphisms of $G$.  In a nonabelian free metabelian nilpotent group $G$
of class $k$, the  group of normal automorphisms is a subgroup of $\textrm{IA}(G)$. 
Since $\textrm{IA}(G)$ stabilizes the lower central series of $G$, it is nilpotent of class $k-1$ (see for instance \cite[p. 9]{SE}), and thus so is the  group of normal automorphisms. Also this result is a consequence of the fact that in a nilpotent group of class $k\geq 2$, the group generated by all  polynomial automorphisms is nilpotent of class $k-1$  \cite{EN3}.\\
If $d>1$ denotes the rank of the free metabelian nilpotent group $G$
of class $k$, note that $\textrm{IA}(G)$ is metabelian if $d=2$ (by a result of C.K. Gupta 
\cite{GU}) or if $k\leq 4$ (in this case, $\textrm{IA}(G)$ is nilpotent of class $\leq 3$).
On the other hand, it is worth pointing out that contrary to the group of normal automorphisms, $\textrm{IA}(G)$ is not metabelian if $d>2$ and $k>4$. 
Indeed, without loss of generality, we can assume that $G$ is the free metabelian nilpotent group of class $5$ and of rank 3, freely generated by $a,b,c$. In this group, we define three IA-automorphisms $f,g,h$ by
$$f(a)=a[a,b],\: f(b)=b,\: f(c)=c,$$
$$g(a)=a,\: g(b)=b[b,c],\: g(c)=c,$$
$$h(a)=a,\: h(b)=b,\: h(c)=c[c,a].$$
It follows
$$[f,g](a)=a[c^{-1},b,a],\: [f,g](b)=b,\: [f,g](c)=c,$$
$$[f,h](a)=a,\: [f,h](b)=b,\: [f,h](c)=c[a,b^{-1},c],$$
whence
$$[f,h]\circ[f,g](c)=c[a,b^{-1},c],$$
$$[f,g]\circ[f,h](c)=c[a,b^{-1},c][c,b,a,b,c].$$
Consequently, $[f,g]$ and $[f,h]$ do not commute, hence  $\textrm{IA}(G)$ is not metabelian.
\section{Proof of Theorem 2.1.}
In all this section, we consider a fixed set $S$ of cardinality $\geq 2$ 
and we denote by $M_{k}$ the free metabelian nilpotent group
of class $k$ freely 
generated by $S$. In other words, $M_{k}=F/F''\gamma_{k+1}(F)$, where 
$F$ is the free group freely generated by $S$ and $\gamma_{k+1}(F)$ 
the $(k+1)$th term of the lower central series of $F$.
The normal closure in a group $G$ of an element $a$ is written
$\langle a^{G}\rangle$.
\begin{lem}
For any distinct elements $a,b\in S$ and 
any integer $t$, the subgroup 
$\langle (a^{t}b)^{M_{k}} \rangle\cap \gamma_{k}(M_{k})\unlhd M_{k}$ is generated 
by the set of elements of the form
$$\left[ a^{t}b,c_{1},\ldots,c_{k-1} \right],\;\; with\; 
c_{1},\ldots,c_{k-1}\in S \;\;(we \; suppose \; k>1). $$
Moreover, for any 
subset $\{ a=a_{0},a_{1},\ldots,a_{r}=b\}\subseteq S$ containing $a$ 
and $b$, the subgroup 
$$\langle (a^{t}b)^{M_{k}} \rangle\cap \gamma_{k}(M_{k})\cap 
\langle a_{0},a_{1},\ldots,a_{r}\rangle$$ 
is generated 
by the set of elements of the form
$\left[ a^{t}b,c_{1},\ldots,c_{k-1} \right]$, with 
$c_{1},\ldots,c_{k-1}\in  \{ a_{0},a_{1},\ldots,a_{r}\}$. 
\end{lem}
\begin{proof} First suppose that $t=0$ and consider an 
element $w\in \langle b^{M_{k}} \rangle\cap 
\gamma_{k}(M_{k})$. Hence $w$ is a product of elements of the form
$[c_{0},c_{1},\ldots,c_{k-1}]^{\pm 1}$, with $c_{i}\in S$. More 
precisely, we can write $w=w_{0}w_{1}$, where $w_{0}$ (resp. $w_{1}$) 
is a product of elements of the form $[c_{0},c_{1},\ldots,c_{k-1}]^{\pm 1}$ 
with $c_{i}\in S\setminus\{ b\}$, (resp. with $c_{i}\in S$, the element $b$ 
occuring once at least in  $[c_{0},c_{1},\ldots,c_{k-1}]$). 
In fact, substituting 1 for the indeterminate $b$ in the relation
$w=w_{0}w_{1}$ and using the fact that $w$ lies in 
$\langle b^{M_{k}} \rangle$, we obtain $w_{0}=1$. 
Thus $w$ is a product of elements 
of the form $[c_{0},c_{1},\ldots,c_{k-1}]^{\pm 1}$, with $c_{i}=b$ for 
some $i\in \{0,\ldots,k-1\}$. If $i=1$, we can write
$[c_{0},b,\ldots,c_{k-1}]=[b,c_{0},\ldots,c_{k-1}]^{-1}$. If $i>1$, 
we have 
$[c_{0},c_{1},\ldots,b,\ldots,c_{k-1}]=[c_{0},c_{1},b,\ldots,c_{k-1}]$;
by using the relation $[c_{0},c_{1},b]=[c_{1},b,c_{0}]^{-1}[b,c_{0},c_{1}]^{-1}$, 
it follows
$$[c_{0},c_{1},\ldots,b,\ldots,c_{k-1}]=
[b,c_{1},c_{0},\ldots,c_{k-1}]
[b,c_{0},c_{1},\ldots,c_{k-1}]^{-1}.$$ 
Thus we have shown that any element of $\langle b^{M_{k}} \rangle\cap 
\gamma_{k}(M_{k})$ is a product of elements of the form 
$\left[ b,c_{1},\ldots,c_{k-1} \right]^{\pm 1}$, with 
$c_{i}\in S$. Since  $\left[ b,c_{1},\ldots,c_{k-1} \right]\in
\langle b^{M_{k}} \rangle\cap \gamma_{k}(M_{k})$, the first 
part of our lemma is proved when $t=0$.\\
Now consider the general case. Actually, since clearly 
$S'=\{ a^tb\}\cup S\setminus\{ b\}$ is a free basis of $M_{k}$, we 
can use the result obtained in the particular case. It follows that 
$\langle (a^{t}b)^{M_{k}} \rangle\cap \gamma_{k}(M_{k})\unlhd M_{k}$ is generated 
by the set of elements of the form
$\left[ a^{t}b,c_{1},\ldots,c_{k-1} \right]$, with 
$c_{i}\in S'$. But in fact we may take $c_{i}\in S$ and so conclude, 
since
$$\left[ a^{t}b,c_{1},\ldots, a^{t}b,\ldots,c_{k-1} \right]=
\left[ a^{t}b,c_{1},\ldots, a,\ldots,c_{k-1} \right]^{t}
\left[ a^{t}b,c_{1},\ldots, b,\ldots,c_{k-1} \right].$$ 
Finally, consider an element $w\in \langle (a^{t}b)^{M_{k}} \rangle\cap \gamma_{k}(M_{k})\cap 
\langle a_{0},a_{1},\ldots,a_{r}\rangle$. 
we can express $w$ in the form $w=w'w''$, where $w'$ (resp. $w''$) 
is a product of elements of the form $[a^tb,c_{1},\ldots,c_{k-1}]^{\pm 1}$ 
with $c_{1},\ldots,c_{k-1}\in \{ a_{0},a_{1},\ldots,a_{r}\}$ 
(resp. with $c_{1},\ldots,c_{k-1}\in S$, an element of
$S\setminus \{ a_{0},a_{1},\ldots,a_{r}\}$ occuring once at least in
the sequence $c_{1},\ldots,c_{k-1}$).
Substituting 1 for all indeterminates in
$S\setminus \{ a_{0},a_{1},\ldots,a_{r}\}$, the equality $w=w'w''$ 
gives $w=w'$. This completes the proof of the lemma. 
\end{proof}
As usual, the expression $[x,_{n}y]$ is defined in a group by 
$[x,_{0}y]=x$ and $[x,_{n}y]=[[x,_{n-1}y],y]$ for each positive 
integer $n$. 

For a fixed subset $\{ a_{0},\ldots,a_{r}\}\subseteq S$ 
and a function $\Delta : \{ 0,\ldots,r\}\to {\mathbb N}$, we define 
in $M_{k}$ the symbol $\left[x,y,\Delta\right]$ ($x,y\in M_{k}$) by  
$$\left[x,y,\Delta\right]=
\left[x,y,_{\Delta(0)}a_{0},_{\Delta(1)}a_{1},\ldots,_{\Delta(r)}a_{r}\right].$$
Note that for any sequence $b_{1},\ldots,b_{k}$ of elements of
$\{ a_{0},\ldots,a_{r}\}$, there is a function 
$\Delta : \{ 0,\ldots,r\}\to {\mathbb N}$ such that 
$\left[x,y,b_{1},\ldots,b_{k}\right]=\left[x,y,\Delta\right]$, with 
$\Delta(0)+\cdots+\Delta(r)=k$ (it suffices to apply the relation
$[u,v,w]=[u,w,v]$, valid in any metabelian group whenever $u$ belongs 
to the derived subgroup). 
If $j,j'$ are distinct given integers in 
$\{ 0,\ldots,r\}$ and if $\Delta(j)\not = 0$, we define the function 
$\Delta_{(j)}^{(j')}: \{ 0,\ldots,r\}\to {\mathbb N}$  
by 
$$\Delta_{(j)}^{(j')}(j)=\Delta(j)-1,\;\;
\Delta_{(j)}^{(j')}(j')=\Delta(j')+1\;\; {\rm and} $$
$$\Delta_{(j)}^{(j')}(i)=\Delta(i)\;\;
{\rm for \; all} \; i\in\{ 0,\ldots,r\}\setminus\{ j,j'\}.$$
When $\Delta$ is not the 
null-function, we shall denote by $m(\Delta)$ 
the least integer $j$ such that $\Delta(j)\not = 0$. 

If $S$ is ordered, we may define in $M_{k}$ {\em basic commutators} 
(see for example \cite[Chapter 3]{NE}. Recall that a basic commutator of 
weight $k'$ ($2\leq k'\leq k$) is a commutator of the form 
$\left[b_{1},b_{2},\ldots,b_{k'}\right]$ ($b_{i}\in S$), with
$b_{1}>b_{2}$ and $b_{2}\leq b_{3}\leq\cdots \leq b_{k'}$. Any set of 
these commutators freely generates a free abelian subgroup of $M_{k}'$.

In the next lemma, we aim to express a product of commutateurs 
of the form $[a_{s},a_{i},\Delta]$ as a product where only basic 
commutators occur.
\begin{lem} 
Let $\{ a_{0},\ldots,a_{r}\}$ be a 
finite subset of $S$ ($r>0$). Choose an integer $s\in \{0,\ldots,r\}$ 
and consider an element $w\in M_{k+2}$ ($k>0$) of the form
$$w=\prod_{i,\Delta}[a_{s},a_{i},\Delta]^{\epsilon(i,\Delta)}\quad 
(\epsilon(i,\Delta)\in {\mathbb Z}),$$
where the product is taken over all integers 
$i\in \{ 0,\ldots,r\}$ and all functions 
$\Delta : \{ 0,\ldots,r\}\to {\mathbb N}$ such that 
$\Delta(0)+\cdots+\Delta(r)=k$.
Then:\\
{\rm (i)} We have
$$w=\prod_{i<s,\, i\leq m(\Delta)}[a_{s},a_{i},\Delta]
^{\epsilon(i,\Delta)}
\prod_{s<i,\, i\leq 
m(\Delta)}[a_{i},a_{s},\Delta]^{-\epsilon(i,\Delta)}$$
$$\times \prod_{s\leq m(\Delta),\, m(\Delta)<i}[a_{i},a_{s},\Delta]
^{-\epsilon(i,\Delta)}
\prod_{m(\Delta)<s,\, m(\Delta)<i}[a_{i},a_{m(\Delta)},\Delta^{(s)}_{(m(\Delta))}]
^{-\epsilon(i,\Delta)}$$
$$ \times \prod_{m(\Delta)<s,\, m(\Delta)<i}
[a_{s},a_{m(\Delta)},\Delta^{(i)}_{(m(\Delta))}]
^{\epsilon(i,\Delta)}$$
(in all these products, $i$ lies in  
$\{0,\ldots,r\}\setminus\{s\}$).\\
{\rm (ii)} We have $w=1$ only if all exponents $\epsilon(i,\Delta)$ 
with $i\in \{0,\ldots,r\}\setminus\{s\}$ occuring in the expression of $w$ are zero.
\end{lem}
\begin{proof} (i) First we write $w$ as a product of two factors:
$$w=\prod_{i\leq m(\Delta)}[a_{s},a_{i},\Delta]^{\epsilon(i,\Delta)}
\prod_{m(\Delta)<i}[a_{s},a_{i},\Delta]^{\epsilon(i,\Delta)}.$$
The first factor can be expressed in the form
\begin{eqnarray*}
\prod_{i\leq m(\Delta)}[a_{s},a_{i},\Delta]^{\epsilon(i,\Delta)}	 & = & 
\prod_{i<s,\, i\leq m(\Delta)}[a_{s},a_{i},\Delta]^{\epsilon(i,\Delta)}  
\prod_{s<i,\, i\leq m(\Delta)}[a_{s},a_{i},\Delta]^{\epsilon(i,\Delta)}\\
	{} & = & \prod_{i<s,\, i\leq m(\Delta)}[a_{s},a_{i},\Delta]^{\epsilon(i,\Delta)}  
\prod_{s<i,\, i\leq m(\Delta)}[a_{i},a_{s},\Delta]^{-\epsilon(i,\Delta)}.
\end{eqnarray*}
Now we deal with the second facteur. We have
\begin{eqnarray*}
\prod_{m(\Delta)<i}[a_{s},a_{i},\Delta]^{\epsilon(i,\Delta)} & = &
\prod_{s\leq m(\Delta),\, m(\Delta)<i}[a_{s},a_{i},\Delta]^{\epsilon(i,\Delta)}
\prod_{m(\Delta)<s,\, m(\Delta)<i}[a_{s},a_{i},\Delta]^{\epsilon(i,\Delta)} \\
{} & = & \prod_{s\leq m(\Delta),\, m(\Delta)<i}[a_{i},a_{s},\Delta]^{-\epsilon(i,\Delta)}
\prod_{m(\Delta)<s,\, m(\Delta)<i}[a_{s},a_{i},\Delta]^{\epsilon(i,\Delta)}. 
\end{eqnarray*}
Therefore Lemma 3.2(i) is proved if we show the relation
\begin{equation}
\prod_{m(\Delta)<s,\, 
m(\Delta)<i}[a_{s},a_{i},\Delta]^{\epsilon(i,\Delta)}=
\end{equation} 
$$ \prod_{m(\Delta)<s,\, m(\Delta)<i}[a_{i},a_{m(\Delta)},\Delta^{(s)}_{(m(\Delta))}]
^{-\epsilon(i,\Delta)}
\times \prod_{m(\Delta)<s,\, m(\Delta)<i}
[a_{s},a_{m(\Delta)},\Delta^{(i)}_{(m(\Delta))}]
^{\epsilon(i,\Delta)}.	
$$
For that, write more explicitly the commutator $[a_{s},a_{i},\Delta]$
(in the following equalities, we write $m$ instead of $m(\Delta)$):
\begin{eqnarray*}
	[a_{s},a_{i},\Delta] & = & 
	\left[a_{s},a_{i},_{\Delta(0)}a_{0},\ldots,_{\Delta(r)}a_{r}\right]  \\
	 & = &  \left[a_{s},a_{i},_{\Delta(m)}a_{m},\ldots,_{\Delta(r)}a_{r}\right] \\
	 & = & \left[a_{s},a_{i},a_{m},_{\Delta(m)-1}a_{m},\ldots,_{\Delta(r)}a_{r}\right].
\end{eqnarray*}
Since 
$\left[a_{s},a_{i},a_{m}\right]=
\left[a_{i},a_{m},a_{s}\right]^{-1}\left[a_{m},a_{s},a_{i}\right]^{-1}
=\left[a_{i},a_{m},a_{s}\right]^{-1}\left[a_{s},a_{m},a_{i}\right]$,
we obtain
\begin{equation}
[a_{s},a_{i},\Delta]=[a_{i},a_{m},\Delta_{(m)}^{(s)}]^{-1}
[a_{s},a_{m},\Delta_{(m)}^{(i)}].	
\end{equation}
Relation (1) is now an immediate consequence of 
(2). \\
(ii) Choose an order in $S$ such that $a_{0}<a_{1}<\cdots <a_{r}$.
Under this condition, all commutators occuring in the expression of 
$w$ obtained in the first part of the lemma are basic commutators.
Suppose that $w=1$. 
Since basic commutators of the form $[a_{j},a_{i},\ldots]$ 
($i,j\in\{0,\ldots,r\}\setminus \{ s\}$) occur only in the fourth 
factor, we have $\epsilon(i,\Delta)=0$ whenever 
$m(\Delta)<s$ and $m(\Delta)<i$. 
It follows 
$$w=\prod_{i<s,\, i\leq m(\Delta)}[a_{s},a_{i},\Delta]
^{\epsilon(i,\Delta)}
\prod_{s<i,\, i\leq 
m(\Delta)}[a_{i},a_{s},\Delta]^{-\epsilon(i,\Delta)}$$
$$\times \prod_{s\leq m(\Delta),\, m(\Delta)<i}[a_{i},a_{s},\Delta]
^{-\epsilon(i,\Delta)}=1.$$
But all basic commutators occuring in this equality are distinct. Thus 
$\epsilon(i,\Delta)=0$ for all integers $i\not = s$ and all functions $\Delta$, as required.
\end{proof}
\begin{lem}
Let $\varphi$ be a normal automorphism 
of $M_{k+2}$ ($k>0$) acting trivially on 
$M_{k+2}/\gamma_{k+2}(M_{k+2})$. Then, for all distinct elements 
$a,b\in S$, there exists a generalized inner automorphism 
$\psi$ of $M_{k+2}$ such that $\varphi(a)=\psi(a)$ and $\varphi(b)=\psi(b)$. 
\end{lem}
\begin{proof} Let $a,b$ be two distinct elements of $S$. Then
$a^{-1}\varphi(a)$ and $b^{-1}\varphi(b)$ belong to 
$\langle a^{M_{k+2}}\rangle \cap \gamma_{k+2}(M_{k+2})$ and
$\langle b^{M_{k+2}}\rangle \cap \gamma_{k+2}(M_{k+2})$ respectively. 
By Lemma 3.1, there is a finite subset
$\{ a=a_{0},a_{1},\ldots,a_{r}=b\}\subseteq S$ such that 
$$
\varphi(a)=\varphi(a_{0})=a_{0}\prod_{i,\Delta}[a_{0},a_{i},\Delta]^{\alpha(i,\Delta)}\quad   
(\alpha(i,\Delta)\in {\mathbb Z}),$$
$$
\varphi(b)=\varphi(a_{r})=a_{r}\prod_{i,\Delta}[a_{r},a_{i},\Delta]^{\beta(i,\Delta)}\quad 
(\beta(i,\Delta)\in {\mathbb Z}),$$
where the two products are taken over all integers 
$i\in \{ 0,\ldots,r\}$ and all functions 
$\Delta : \{ 0,\ldots,r\}\to {\mathbb N}$ with 
$\Delta(0)+\cdots+\Delta(r)=k$ (as in Lemma 3.2,  $[a_{s},a_{i},\Delta ]$ is 
defined by
$[a_{s},a_{i},\Delta]=[a_{s},a_{i},_{\Delta(0)}a_{0},\ldots,_{\Delta(r)}a_{r}]$). 
Note that if $|S|=2$ (and so $r=1$), Lemma 3.3 is easily verified by taking the 
generalized inner automorphism $\psi$ defined by
$$\psi(x)=x\prod_{\Delta}[x,a_{1},\Delta]^{\alpha(1,\Delta)} 
\prod_{\Delta}[x,a_{0},\Delta]^{\beta(0,\Delta)}.$$
Thus we can assume in the following that $|S|>2$.
For any positive integer $t$, $(a_{0}^ta_{r})^{-1}\varphi(a_{0}^ta_{r})=
(a_{0}^ta_{r})^{-1}\varphi(a_{0})^t\varphi(a_{r})$ belongs to 
$\langle a_{0},a_{1},\ldots,a_{r}\rangle$. Therefore, once again by Lemma 
2.1, $\varphi(a_{0}^ta_{r})$ can be expressed in the form
$$\varphi(a_{0}^ta_{r})=a_{0}^ta_{r}\prod_{i,\Delta}[a_{0}^ta_{r},a_{i},\Delta]^{\eta_{t}(i,\Delta)}\quad   
(\eta_{t}(i,\Delta)\in {\mathbb Z}).$$
Since $[a_{0},a_{i},\Delta]$ is a commutator of weight $k+2$ in a 
nilpotent group of class $k+2$, we obtain
$$\varphi(a_{0}^ta_{r})=a_{0}^ta_{r}\prod_{i,\Delta}[a_{0},a_{i},\Delta]^{t\eta_{t}(i,\Delta)}
\prod_{i,\Delta}[a_{r},a_{i},\Delta]^{\eta_{t}(i,\Delta)}.$$
Thus the relation $\varphi(a_{0}^ta_{r})=\varphi(a_{0})^t\varphi(a_{r})$ implies that 
\begin{equation}
\prod_{i,\Delta}[a_{0},a_{i},\Delta]^{t\eta_{t}(i,\Delta)}
\prod_{i,\Delta}[a_{r},a_{i},\Delta]^{\eta_{t}(i,\Delta)}
\end{equation} 
$$=\, \prod_{i,\Delta}[a_{0},a_{i},\Delta]^{t\alpha(i,\Delta)}
\prod_{i,\Delta}[a_{r},a_{i},\Delta]^{\beta(i,\Delta)}. $$
Choose an order in $S$ such that $a_{0}<a_{1}<\cdots <a_{r}$.
Then we can use Lemma 3.2(i) (with $s=0$ or $s=r$) to express each product in (3) as a 
product of basic commutators (or their inverses). The first product gives
$$\prod_{i,\Delta}[a_{0},a_{i},\Delta]^{t\eta_{t}(i,\Delta)}
=\prod_{0<i,\, i\leq m(\Delta)}[a_{i},a_{0},\Delta]^{-t\eta_{t}(i,\Delta)}
\prod_{0\leq m(\Delta),\, m(\Delta)<i}
[a_{i},a_{0},\Delta]^{-t\eta_{t}(i,\Delta)}$$
and so
$\displaystyle \prod_{i,\Delta}[a_{0},a_{i},\Delta]^{t\eta_{t}(i,\Delta)}
=\prod_{i\not = 0,\, \Delta}
[a_{i},a_{0},\Delta]^{-t\eta_{t}(i,\Delta)},$
where this product is taken over all integers 
$i\in \{ 1,\ldots,r\}$ and all functions 
$\Delta : \{ 0,\ldots,r\}\to {\mathbb N}$ with 
$\Delta(0)+\cdots+\Delta(r)=k$.\\
Likewise we have $\displaystyle 
\prod_{i,\Delta}[a_{0},a_{i},\Delta]^{t\alpha(i,\Delta)}
=\prod_{i\not = 0,\, \Delta}
[a_{i},a_{0},\Delta]^{-t\alpha(i,\Delta)}$.\\
In the same way, applying Lemma 3.2(i) with $s=r$,
the second product of (3) gives
$$\prod_{i,\Delta}[a_{r},a_{i},\Delta]^{\eta_{t}(i,\Delta)}=
\prod_{i\not = r,\, i\leq m(\Delta)}[a_{r},a_{i},\Delta]
^{\eta_{t}(i,\Delta)}$$
$$\times \prod_{i\not = r,\, m(\Delta)<i}[a_{i},a_{m(\Delta)},\Delta^{(r)}_{(m(\Delta))}]
^{-\eta_{t}(i,\Delta)}
\prod_{i\not = r,\, m(\Delta)<i}
[a_{r},a_{m(\Delta)},\Delta^{(i)}_{(m(\Delta))}]
^{\eta_{t}(i,\Delta)}.$$ 
Finally, the last product of (3) gives 
$$\prod_{i,\Delta}[a_{r},a_{i},\Delta]^{\beta(i,\Delta)}=
\prod_{i\not = r,\, i\leq m(\Delta)}[a_{r},a_{i},\Delta]
^{\beta(i,\Delta)}$$
$$\times \prod_{i\not = r,\, m(\Delta)<i}[a_{i},a_{m(\Delta)},\Delta^{(r)}_{(m(\Delta))}]
^{-\beta(i,\Delta)}
\prod_{i\not = r,\, m(\Delta)<i}
[a_{r},a_{m(\Delta)},\Delta^{(i)}_{(m(\Delta))}]
^{\beta(i,\Delta)}.$$
Thus relation (3) can be written in the form
\begin{equation}
\prod_{i\not = 0,\, \Delta}
[a_{i},a_{0},\Delta]^{-t\eta_{t}(i,\Delta)}
\prod_{i\not = r,\, i\leq m(\Delta)}[a_{r},a_{i},\Delta]
^{\eta_{t}(i,\Delta)}
\end{equation}
$$\times \prod_{i\not = r,\, m(\Delta)<i}[a_{i},a_{m(\Delta)},\Delta^{(r)}_{(m(\Delta))}]
^{-\eta_{t}(i,\Delta)}
\prod_{i\not = r,\, m(\Delta)<i}
[a_{r},a_{m(\Delta)},\Delta^{(i)}_{(m(\Delta))}]
^{\eta_{t}(i,\Delta)}$$ 
$$=\quad \prod_{i\not = 0,\, \Delta}
[a_{i},a_{0},\Delta]^{-t\alpha(i,\Delta)}
\prod_{i\not = r,\, i\leq m(\Delta)}[a_{r},a_{i},\Delta]
^{\beta(i,\Delta)}$$
$$
\times \prod_{i\not = r,\, m(\Delta)<i}[a_{i},a_{m(\Delta)},\Delta^{(r)}_{(m(\Delta))}]
^{-\beta(i,\Delta)}
\prod_{i\not = r,\, m(\Delta)<i}
[a_{r},a_{m(\Delta)},\Delta^{(i)}_{(m(\Delta))}]
^{\beta(i,\Delta)}	
$$
(so each commutator occuring in this relation is a basic commutator).\\
Now consider an integer $i\in \{ 1,\ldots,r-1\}$ and a function
$\Delta : \{ 0,\ldots,r\}\to {\mathbb N}$, with 
$\Delta(0)+\cdots+\Delta(r)=k$ (we can always suppose that $r>1$ since 
$|S|>2$).
By identifying the exponents of the basic commutator
$[a_{i},a_{0},\Delta]$ of each side of relation (4), it is easy to see 
that
\begin{equation}
t\eta_{t}(i,\Delta)+\eta_{t}\left(i,\Delta_{(r)}^{(0)}\right)=
t\alpha(i,\Delta)+\beta\left(i,\Delta_{(r)}^{(0)}\right)	
\end{equation} 
if $\Delta(r)>0$, and 
$\eta_{t}(i,\Delta)=\alpha(i,\Delta)$ if $\Delta(r)=0$. 
We prove by induction on $\Delta(r)$ that actually, we have always 
the equality $\eta_{t}(i,\Delta)=\alpha(i,\Delta)$. 
At first observe that if 
$\Delta(r)>0$, we have $\Delta_{(r)}^{(0)}(r)=\Delta(r)-1$ and so 
$\eta_{t}\left(i,\Delta_{(r)}^{(0)}\right)=\alpha\left(i,\Delta_{(r)}^{(0)}\right)$ by 
induction. Hence relation (5) implies that
$$\alpha\left(i,\Delta_{(r)}^{(0)}\right)-\beta\left(i,\Delta_{(r)}^{(0)}\right)=
t\left\{\alpha(i,\Delta)-\eta_{t}(i,\Delta)\right\}.$$
Consequently, each positive integer $t$ divides the integer
$\alpha\left(i,\Delta_{(r)}^{(0)}\right)-\beta\left(i,\Delta_{(r)}^{(0)}\right)$,
which is independant of $t$.
It follows that 
$\alpha\left(i,\Delta_{(r)}^{(0)}\right)=\beta\left(i,\Delta_{(r)}^{(0)}\right)$ and so 
$\eta_{t}(i,\Delta)=\alpha(i,\Delta)$, as required.\\
Using theses relations and taking $t=1$, relation (3) implies
$$\prod_{\Delta}[a_{0},a_{r},\Delta]^{\eta(r,\Delta)}
\prod_{i,\Delta}[a_{r},a_{i},\Delta]^{\eta(i,\Delta)}
=\prod_{\Delta}[a_{0},a_{r},\Delta]^{\alpha(r,\Delta)}
\prod_{i,\Delta}[a_{r},a_{i},\Delta]^{\beta(i,\Delta)}$$
(we write $\eta$ for $\eta_{1}$) and so 
$$\prod_{i,\Delta}[a_{r},a_{i},\Delta]^{\beta(i,\Delta)}=
\prod_{i,\Delta}[a_{r},a_{i},\Delta]^{\eta(i,\Delta)}
\prod_{\Delta}[a_{r},a_{0},\Delta]^{\alpha(r,\Delta)-\eta(r,\Delta)}.$$
Since 
$\displaystyle \varphi(a_{r})=
a_{r}\prod_{i,\Delta}[a_{r},a_{i},\Delta]^{\beta(i,\Delta)}$, 
we obtain
\begin{equation}
\varphi(a_{r})=a_{r}\prod_{i,\Delta}[a_{r},a_{i},\Delta]^{\eta(i,\Delta)}
\prod_{\Delta}[a_{r},a_{0},\Delta]^{\alpha(r,\Delta)-\eta(r,\Delta)}	
\end{equation}
Now consider the generalized inner automorphism $\psi$ defined by
$$\psi(x)=x\prod_{i=1,\ldots,r,\, \Delta}
[x,a_{i},\Delta]^{\alpha(i,\Delta)}
\prod_{\Delta}[x,a_{0},\Delta]^{\alpha(r,\Delta)+\eta(0,\Delta)-\eta(r,\Delta)}.$$
We have
$$\psi(a_{0})=a_{0}\prod_{i=1,\ldots,r,\, \Delta}
[a_{0},a_{i},\Delta]^{\alpha(i,\Delta)}=\varphi(a_{0}).$$
In the same way,
\begin{eqnarray*}
	\psi(a_{r}) & = & a_{r}\prod_{i=1,\ldots,r-1,\, \Delta}
[a_{r},a_{i},\Delta]^{\alpha(i,\Delta)}
\prod_{\Delta}[a_{r},a_{0},\Delta]^{\alpha(r,\Delta)+\eta(0,\Delta)-\eta(r,\Delta)}  \\
	{} & = & a_{r}\prod_{i=1,\ldots,r-1,\, \Delta}
[a_{r},a_{i},\Delta]^{\eta(i,\Delta)}
\prod_{\Delta}[a_{r},a_{0},\Delta]^{\alpha(r,\Delta)+\eta(0,\Delta)-\eta(r,\Delta)}   \\
	{} & = & a_{r}\prod_{i,\, \Delta}
[a_{r},a_{i},\Delta]^{\eta(i,\Delta)}
\prod_{\Delta}[a_{r},a_{0},\Delta]^{\alpha(r,\Delta)-\eta(r,\Delta)} 
\end{eqnarray*}
and so $\psi(a_{r})=\varphi(a_{r})$ by (6). This completes the proof of 
Lemma 3.3. 
\end{proof}
Lemma 3.3 is actually a weak form of the next result.
\begin{lem}
Let $\varphi$ be a normal automorphism 
of $M_{k+2}$ ($k>0$) acting trivially on 
$M_{k+2}/\gamma_{k+2}(M_{k+2})$. Then $\varphi$ is
a generalized inner automorphism of $M_{k+2}$.
\end{lem}
\begin{proof} We can assume that $|S|>2$ (otherwise Lemma 3.4
is a consequence of Lemma 3.3). Consider two distinct elements 
$a,b\in S$. According to Lemma 3.3, there exists a generalized inner 
automorphism $\psi$ such that $\varphi(a)=\psi(a)$ and $\varphi(b)=\psi(b)$.
It suffices to prove that for any $c\in S\setminus\{ a,b\}$, we have $\varphi(c)=\psi(c)$.
For that, apply again Lemma 3.3: there are generalized inner 
automorphisms $\psi',\psi''$ such that 
$\varphi(a)=\psi'(a)$, $\varphi(c)=\psi'(c)$ and $\varphi(b)=\psi''(b)$,
$\varphi(c)=\psi''(c)$. There exists a finite subset
$\{ a_{0},\ldots,a_{r}\}\subseteq S$, containing $a,b,c$, such that
$\psi,\psi',\psi''$ can be defined by the equations
\begin{eqnarray*}
	\psi(x) & = & 
	x\prod_{i,\,\Delta}[x,a_{i},\Delta]^{\epsilon(i,\Delta)}\\
	\psi'(x) & = & 
	x\prod_{i,\,\Delta}[x,a_{i},\Delta]^{\epsilon'(i,\Delta)} \\
	\psi''(x) & = & x\prod_{i,\,\Delta}[x,a_{i},\Delta]^{\epsilon''(i,\Delta)}
\end{eqnarray*}
(the products are taken over all integers 
$i\in \{ 0,\ldots,r\}$ and all functions 
$\Delta : \{ 0,\ldots,r\}\to {\mathbb N}$ with 
$\Delta(0)+\cdots+\Delta(r)=k$). Since $\psi(a)=\psi'(a)$, we have 
$$a\prod_{i,\,\Delta}[a,a_{i},\Delta]^{\epsilon(i,\Delta)}=
a\prod_{i,\,\Delta}[a,a_{i},\Delta]^{\epsilon'(i,\Delta)}$$ and so
$$\prod_{i,\,\Delta}[a,a_{i},\Delta]^{\epsilon(i,\Delta)-\epsilon'(i,\Delta)}=1.$$
Applying Lemma 3.2(ii), we obtain 
$\epsilon(i,\Delta)=\epsilon'(i,\Delta)$ for all functions $\Delta$ and 
all integers $i\in \{0,\ldots,r\}$ such that $a_{i}\not = a$.
Similarly, we have $\epsilon(i,\Delta)=\epsilon''(i,\Delta)$ if  
$a_{i}\not = b$ and $\epsilon'(i,\Delta)=\epsilon''(i,\Delta)$ if  
$a_{i}\not = c$. It follows that 
$\epsilon(i,\Delta)=\epsilon'(i,\Delta)$ for all function $\Delta$ and 
all integer $i\in \{0,\ldots,r\}$, hence $\psi=\psi'$.
Thus $\varphi(c)=\psi'(c)=\psi(c)$, as required. 
\end{proof}
\begin{proof}[Proof of Theorem 2.1] We argue by induction on the 
nilpotency class $k$ of $M_{k}$. If $k=2$, the result follows from 
\cite[Theorem 2(ii)]{EN2} (in this case, each normal automorphism is
inner). Now consider a normal automorphism $\varphi$ of $M_{k}$, with 
$k>2$. Then $\varphi$ induces a normal automorphism on the quotient 
group $M_{k}/\gamma_{k}(M_{k})$. By induction, since this quotient is isomorphic to 
$M_{k-1}$, there exists a generalized inner automorphism $\psi:M_{k}\to M_{k}$
such that $\varphi(x)=\psi(x)\theta(x)$, where $\theta(x)$ is an 
element of $\gamma_{k}(M_{k})$. It follows that 
$\psi^{-1}(\varphi(x))=x\psi^{-1}(\theta(x))$. Thus  
$\psi':=\psi^{-1}\circ\varphi$ is a normal automorphism 
of $M_{k}$ acting trivially on 
$M_{k}/\gamma_{k}(M_{k})$. By Lemma 3.4, $\psi'$ is a
generalized inner automorphism, and so is $\varphi=\psi\circ\psi'$. 
This completes the proof of Theorem 2.1.
\end{proof}
We end with a question: assume $\varphi$ is a normal automorphism of 
a free nilpotent group $G$ of class $k\geq 4$ (and of rank $>1$). Is $\varphi$ of 
the form
$$\varphi:x\to 
x\left[x,u_{1}\right]^{\lambda(1)}\ldots\left[x,u_{m}\right]^{\lambda(m)}\;\;\;
(u_{i}\in G,\; \lambda(i)\in{\mathbb Z})\, ?$$

%
%
  
%
\end{document}